\documentclass[12pt,a4paper]{amsart}
\usepackage{amssymb,latexsym,verbatim}
\usepackage{mathrsfs}

\headsep=1cm \numberwithin{equation}{section}
\newtheorem{thm}{Theorem}[section]
 \newtheorem{cor}[thm]{Corollary}
 \newtheorem{lem}[thm]{Lemma}
 \newtheorem{prop}[thm]{Proposition}
 
 \newtheorem{defn}[thm]{Definition}
 \newtheorem{exam}[thm]{Example}
 
 \theoremstyle{definition}
 \theoremstyle{remark}

 \numberwithin{equation}{section}
 \newcommand{\p}{\frak p}
\newcommand\Supp{\operatorname{Supp}}
\newcommand\Ass{\operatorname{Ass}}

\newcommand\Spec{\operatorname{Spec}}

\newcommand\Hom{\operatorname{Hom}}
\newcommand\Ext{\operatorname{Ext}}

\begin{document}
\title[ Cofiniteness of weakly Laskerian local cohomology modules  ]
{Cofiniteness of weakly Laskerian local cohomology modules}
\author[Moharram Aghapournahr and Kamal Bahmanpour]{Moharram Aghapournahr$^*$ and Kamal Bahmanpour\\\\\\\,
\vspace*{0.5cm}Dedicated to Professor Leif Melkersson}

\address{ Department of Mathematic, Faculty of Science,
 Arak University, Arak, 38156-8-8349,
 Iran.} \email{m-aghapour@araku.ac.ir}

\address{Department of Mathematics, Islamic Azad University-Ardabil branch, P.O. Box 5614633167,
Ardabil, Iran.} \email{bahmanpour.k@gmail.com}

\thanks{ 
$^*$Corresponding author: e-mail: {\it m-aghapour@araku.ac.ir} (Moharram Aghapournahr)}%
%\subjclass{}%

%\date{}%
%\dedicatory{}
%\commby{}%
% ----------------------------------------------------------------
\begin{abstract}
Let  $I$ be an ideal of a Noetherian ring R and M be a finitely generated
R-module. We introduce  the class of extension modules of finitely generated
modules by the class of all modules $T$  with $\dim T\leq n$  and we show it by
${\rm FD_{\leq n}}$ where $n\geq -1$ is an integer. We prove that for any  ${\rm
FD_{\leq 0}}$(or minimax) submodule N of $H^t_I(M)$ the R-modules ${\rm
Hom}_R(R/I,H^{t}_I(M)/N)\,\,\, {\rm and}\,\,\,{\rm Ext}^1_R(R/I,H^{t}_I(M)/N)$ are
finitely generated, whenever the modules $H^0_I(M)$, $H^1_I(M)$, ...,
$H^{t-1}_I(M)$ are ${\rm FD_{\leq 1}}$  ( or weakly Laskerian).  As a
consequence, it follows that the associated primes of $H^{t}_I(M)/N$  are finite.
This generalizes  the main results of Bahmanpour and Naghipour \cite{BN} and
\cite{BN2}, Brodmann and Lashgari  \cite{BL}, Khashyarmanesh and Salarian
\cite{KS} and Hong Quy \cite{hung}. We also show that the category  $\mathscr
{FD}^1(R,I)_{cof}$ of $I$-cofinite ${\rm FD_{\leq1}}$ ~ $R$-modules forms an
Abelian subcategory of the category of all $R$-modules.
\end{abstract}
\maketitle
% ----------------------------------------------------------------
\section{Introduction}

The following conjecture was made by Grothendieck in \cite{SGA2}:\\

\noindent {\bf Conjecture:}  \emph{For any ideal $I$ of a Noetherian ring $R$ and
any finite $R$--module $M$, the module $\Hom_{R}(R/I,H^j_I(M))$ is finitely
generated for all $j\ge 0$}.\\

 Here, $H^j_I(M)$
denotes the $j^{th}$  local cohomology module of $M$ with support in $I$.
 Although the conjecture is not true in general as  was shown by
Hartshorne in \cite{Ha}, there are some attempts to show that under some
conditions, for
  some number $t$, the module
$\Hom_{R}(R/I,H^t_I(M))$ is finite, see \cite[Theorem 3.3]{AKS}, \cite[Theorem
6.3.9]{DY2}, \cite[Theorem 2.1]{DY}, \cite[Theorem 2.6]{BN} and \cite[Theorem
2.3]{BN2}.

 In \cite{Ha}, Hartshorne  defined an $R$-module $L$ to
be $I$-{\it cofinite}, if ${\rm Supp}(L)\subseteq V(I)$ and ${\rm Ext}^{i}_{R}(R/I, L)$
is finitely generated module for all $i$. He asked:\\

\emph{If $I$ is an ideal of $R$ and $M$ is a finitely generated $R$-module,
when is $H^t_I(M)$  $I$--cofinite for all i ?}\\

In this direction in section 3 we generalize  \cite[Theorem 3.3]{AKS},
\cite[Theorem 2.6]{BN} and \cite[Theorem 2.3]{BN2} to the class of extension
modules of finitely generated modules by the class of all modules $T$  with
$\dim T\leq 1$ (${\rm FD_{\leq 1}}$). Note that the class of weakly Laskerian
modules is contained in the class of  ${\rm FD_{\leq 1}}$ modules. More
precisely, we shall show that:

\begin{thm}
 Let $R$ be a Noetherian ring and $I$ an ideal of $R$. Let $M$ be
a finitely generated $R$-module and $t\geq 1$ be a positive integer such that the
$R$-modules $H^i_I(M)$ are ${\rm FD_{\leq1}}$ ~ $R$-modules {\rm{(}}or weakly
Laskerian{\rm{)}} for all $i<t$. Then, the following conditions hold:

{\rm(i)} The $R$-modules $H^i_I(M)$ are $I$-cofinite for all $i<t$.

 {\rm(ii)} For all ${\rm FD_{\leq 0}}$ {\rm{(}}or minimax{\rm{)}} submodule $N$ of
$H^{t}_I(M)$, the $R$-modules $${\rm Hom}_R(R/I,H^{t}_I(M)/N)\,\,\, {\rm
and}\,\,\,{\rm Ext}^1_R(R/I,H^{t}_I(M)/N)$$ are finitely generated.
 \end{thm}

As an immediate consequence we prove the following corollary that is a
generalization of Bahmanpour-Naghipour's results in \cite{BN}  and also the
Delfno-Marley's result in \cite{DM} and Yoshida's result in \cite{Yo} for an
arbitrary Noetherian ring.

\begin{cor}
\label{cof1} Let $R$ be a Noetherian ring and $I$ an ideal of $R$. Let $M$ be a
finitely generated $R$-module  such that the $R$-modules $H^i_I(M)$ are ${\rm
FD_{\leq1}}$ {\rm{(}}or weakly Laskerian{\rm{)}}~ $R$-modules for all $i$. Then,

{\rm(i)} the $R$-modules $H^i_I(M)$ are $I$-cofinite for all $i$.

 {\rm(ii)} For any $i\geq 0$ and for any   ${\rm FD_{\leq 0}}$ {\rm{(}}or minimax{\rm{)}} submodule $N$ of
$H^{i}_I(M)$, the $R$-module $H^{i}_I(M)/N$ is  $I$-cofinite.
\end{cor}

 Abazari and Bahmanpour in \cite{AB}  studied cofiniteness of
extension functors of cofinite modules as a generalization of Huneke-Koh's
results in \cite{HK}. In Corollary \ref{torext} we generalis the results of Abazari
and Bahmanpour.

  Hartshorn also posed the following question:\\

{\it Whether the category $\mathscr{M}(R, I)_{cof}$ of $I$-cofinite modules forms
an Abelian subcategory of the category of all $R$-modules? That is, if $f:
M\longrightarrow N$ is an $R$-module homomorphism of $I$-cofinite modules,
are $\ker f$ and ${\rm coker}
f$ $I$-cofinite?}\\

Hartshorne proved that if $I$ is a prime ideal of dimension one in a complete
regular local ring $R$, then the answer to his question is yes. On the other hand,
in \cite{DM}, Delfino and Marley extended this result to arbitrary complete local
rings. Recently, Kawasaki \cite{Ka2} generalized the Delfino and Marley's result
for an arbitrary ideal $I$ of dimension one in a local ring $R$. Finally, more
recently, Sedghi, Bahmanpour and Naghipour in \cite{BNS}  completely have
removed local assumption on $R$. One of the main results of this section is to
prove that the class of  $I$-{\it cofinite} ~${\rm FD_{\leq1}}$  modules compose
an Abelian category (see Theorem \ref{abel}).

%In section 4 we will prove similar results of section 3 concerning the class of
%${\rm FD_{\leq2}}$ modules when $R$ is a local ring.

Let $R$ denote a commutative Noetherian ring, and  let $I$ be an ideal of $R$.
Throughout this paper, $R$ will always be a commutative Noetherian ring with
non-zero identity and $I$ will be an ideal of $R$.  We denote $\{\frak p \in {\rm
Spec}\,R:\, \frak p\supseteq \frak a \}$ by $V(\frak a)$. For any unexplained
notation and terminology we refer the reader to \cite{BH} and \cite{Mat}.

\section{Preliminaries}

%------------------------------------------------------------------------------------------------------------------

Yoshizawa in \cite[Definition 2.1]{Y} defined classes of extension modules of
Serre subcategory by another one as below.

\begin{defn}
\label{D:FD1}
 Let $\mathcal S_1$ and  $\mathcal S_2$ be Serre subcategories of the category of all R-modules. We denote by
 $(\mathcal S_1, \mathcal S_2)$ the class of all $R$-modules M with some $R$-modules
 $S_1\in \mathcal S_1$ and  $S_2 \in \mathcal S_2$
 such that a sequence $0\longrightarrow S_1\longrightarrow  M \longrightarrow S_2
\longrightarrow 0$ is exact.
\end{defn}
We will denote the class of all modules $M$ with $\dim M\leq n$ by  ${\rm
D_{\leq n}}$ and the class of extension modules of finitely generated modules by
the class of ${\rm D_{\leq n}}$ modules by  ${\rm FD_{\leq n}}$ where $n\geq -1$
is an integer. Note that the class of ${\rm FD_{\leq -1}}$ is the same as finitely
generated $R$-modules.
 Recall that a module $M$ is a \emph{minimax} module if there is a
 finitely generated submodule $N$ of $M$ such that the quotient module $M/N$
is artinian. Thus the class of minimax modules is the class of extension modules
of finitely generated modules by the class of Artinian modules. Minimax modules
have been studied by Zink in \cite{Zi} and Z\"{o}schinger in
 \cite{Zrmm,Zrrad}. See also \cite{Ru}.  Recall too that an
$R$-module $M$ is called \emph{weakly Laskerian} if ${\rm Ass}(M/N)$ is a
finite set for each  submodule $N$ of $M$. The class of weakly Laskerian
modules introduced in \cite{DiM}, by Divaani-Aazar and Mafi. Recently, Hung
Quy \cite{hung}, introduced the class of extension modules of finitely generated
modules by the class of all modules  of finite support and named it ${\rm FSF}$
modules. By the following theorem over a Noetherian ring $R$ an $R$-module
$M$ is weakly Laskerian if and
only if is ${\rm FSF}$. \\

\begin{thm}
\label{fsf} Let $R$ be a Noetherian ring and $M$ a nonzero $R$-module. The
following statements are equivalent:

\begin{enumerate}
  \item $M$ is a weakly Laskerian module;
  \item M is an ${\rm FSF}$ module.
\end{enumerate}
\end{thm}
\proof See \cite[Theorem 3.3]{B}.\qed\\

\begin{lem}
\label{FDn}Let $R$ be a Noetherian ring.  Then the  following conditions hold:

{\rm(i)}  Any finitely generated $R$-module and any ${\rm D_{\leq n}}$
$R$-module are  ${\rm FD_{\leq n}}$.

{\rm(ii)} The class of  ${\rm FD_{\leq n-1}}$ modules is contained in the class of
${\rm FD_{\leq n}}$  modules for all $n\geq 0$.

{\rm(iii)} The class of minimax modules is contained in the class of  ${\rm
FD_{\leq 0}}$ that is the class of extension modules of finitely generated
modules by semiartinian modules.

{\rm(iv)} The class of weakly Laskerian modules is contained in the class of ${\rm
FD_{\leq 1}}$.

{\rm(v)} The class of  ${\rm FD_{\leq n}}$~$R$-modules forms a Serre
subcategory of the category of all $R$-modules.
 \end{lem}
\proof (i),   (ii),  (iii) are trivial.

 (iv) Use Theorem \ref{fsf}.

 (v) See \cite[Corollary 4.3 or 4.5]{Y}.\qed\\

\begin{exam}
{\rm{(}i\rm{)}} Let $R$ be a Notherian ring with $\dim R \geq 2$ and let $\p\in
\Spec(R)$ such that
 $\dim R/\p =1$. Let $M=R\oplus E( R/\p)$. It is easy to see that $M$ is an  ${\rm FD_{\leq 1}}$
 $R$-module  that is neither finitely generated nor  ${\rm
D_{\leq 1}}$.\\
{\rm{(}ii\rm{)}} Suppose the set $\Omega$ of maximal ideals of $R$ is infinite.
Then the module $\oplus_{\frak m\in\Omega}R/\frak m$  is ${\rm FD_{\leq 0}}$
module and thus
 ${\rm FD_{\leq 1}}$ but it is not a weakly Laskerian module.

 \end{exam}

\begin{prop}
\label{D1}
 Let $I$ be an ideal of a Noetherian ring $R$ and $M$   be a ${\rm
D_{\leq 1}}$ module such that  $\Supp M\subseteq V(I)$. Then the following
statements are equivalent:

{\rm(i)} $M$ is $I$-cofinite,

{\rm(ii)} The $R$-modules $\Hom_R(R/I,M)$ and $\Ext^1_R(R/I,M)$ are finitely
generated.
 \end{prop}
\proof See \cite[Proposition 2.6]{BNS}.\qed\\

\section{Cofinitness of local cohomology}

In what follows the next theorem plays an important role.

\begin{thm}
\label{FD1}
 Let $I$ be an ideal of a Noetherian ring $R$ and $M$ be an  ${\rm FD_{\leq1}}$  $R$-module
  such that $\Supp M\subseteq V(I)$. Then the following statements
are equivalent:

{\rm(i)} $M$ is $I$-cofinite,

{\rm(ii)} The $R$-modules $\Hom_R(R/I,M)$ and $\Ext^1_R(R/I,M)$ are finitely
generated.
\end{thm}
\proof $(i)\Rightarrow (ii)$ is clear. In order to prove $(ii)\Rightarrow (i)$, by
Definition there is a finitely generated submodule $N$ of $M$ such that the
$R$-module  ${\rm dim}(M/N)\leq 1$ and $\Supp M/N \subseteq V(I)$. Also, the
exact sequence
$$0\rightarrow N \rightarrow M \rightarrow M/N \rightarrow 0,\,\,\,\,\,\,(*)$$
induces the following exact sequence
$$0\longrightarrow {\rm Hom}_R(R/I,N)\longrightarrow{\rm
Hom}_R(R/I,M)\longrightarrow {\rm Hom}_R(R/I,M/N)$$$$\longrightarrow {\rm
Ext}^1_R(R/I,N)\longrightarrow{\rm Ext}^1_R(R/I,M)\longrightarrow {\rm
Ext}^1_R(R/I,M/N)\longrightarrow {\rm Ext}^2_R(R/I,N).$$ Whence, it follows that
the $R$-modules $\Hom_R(R/I,M/N)$ and $\Ext^1_R(R/I,M/N)$ are finitely
generated. Therefore, in view of Proposition \ref{D1}, the $R$-module $M/N$ is
$I$-cofinite. Now it follows from the exact sequence $(*)$ that $M$ is $I$-cofinite.\qed\\

%\begin{cor}
%\label{wl}
% Let $I$ be an ideal of a Noetherian ring $R$ and $M$ be a weakly
%Laskerian $R$-module such that $\Supp M\subseteq V(I)$. Then the following
%statements are equivalent:

%{\rm(i)} $M$ is $I$-cofinite,

%{\rm(ii)} the $R$-modules $\Hom_R(R/I,M)$ and $\Ext^1_R(R/I,M)$ are
%finitely generated.
%\end{cor}
%\proof $(i)\Rightarrow (ii)$ is clear. In order to prove
%$(ii)\Rightarrow (i)$, by Theorem 3.3 there is a finitely generated
%submodule $N$ of $M$ such that the $R$-module $M/N$ has finite
%support. So ${\rm dim}(M/N)\leq 1$ and $\Supp M/N \subseteq V(I)$.
%Also, the exact sequence
%$$0\rightarrow N \rightarrow M \rightarrow M/N \rightarrow 0,\,\,\,\,\,\,(*)$$
%induces the following exact sequence
%$$0\longrightarrow {\rm Hom}_R(R/I,N)\longrightarrow{\rm
%Hom}_R(R/I,M)\longrightarrow {\rm Hom}_R(R/I,M/N)$$$$\longrightarrow {\rm
%Ext}^1_R(R/I,N)\longrightarrow{\rm Ext}^1_R(R/I,M)\longrightarrow {\rm
%Ext}^1_R(R/I,M/N)\longrightarrow {\rm Ext}^2_R(R/I,N).$$ Whence, it follows that
%the $R$-modules $\Hom_R(R/I,M/N)$ and $\Ext^1_R(R/I,M/N)$ are finitely
%generated. Therefore, in view of Proposition 4.1, the $R$-module $M/N$ is
%$I$-cofinite. Now it follows from the exact sequence $(*)$ that $M$ is $I$-cofinite.\qed\\

\begin{lem}
 \label{2.2}
Let $I$ be an ideal of Noetherian ring $R$, $M$ a non-zero
$R$-module and $t\in\Bbb{N}_0$. Suppose that the $R$-module
$H^i_I(M)$ is $I$-cofinite for all $i=0,...,t-1$, and the
$R$-modules ${\rm Ext}^{t}_R(R/I,M)$ and ${\rm Ext}^{t+1}_R(R/I,M)$
are finitely generated. Then the $R$-modules ${\rm
Hom}_R(R/I,H^{t}_I(M))$ and ${\rm Ext}^{1}_R(R/I,H^{t}_I(M))$ are
finitely generated.
 \end{lem}
\proof
See \cite[Theorem 2.1]{DY}  and  \cite[Theorem A]{DY1}.\qed\\

\begin{lem}
\label{FD0}
 Let $I$ be an ideal of a Noetherian ring $R$ and $M$ be an  ${\rm FD_{\leq0}}$  $R$-module
  such that $\Supp M\subseteq V(I)$. Then the following statements
are equivalent:

{\rm(i)} $M$ is $I$-cofinite,

{\rm(ii)} The $R$-module $\Hom_R(R/I,M)$ is finitely generated.
\end{lem}
\proof The proof is similar to the proof of \cite[Proposition 4.3]{Mel}.\qed\\

We are now ready to state and prove the following main results (Theorem
\ref{homext} and   the Corollaries \ref{cof} and \ref{wl}) which are extension of
Bahmanpour-Naghipour's results in \cite{BN} and \cite{BN2},
Brodmann-Lashgari's result in  \cite{BL}, Khashyarmanesh-Salarian's result in
\cite{KS}, Hong Quy's result in  \cite{hung}, and also the Delfno-Marley's result in
\cite{DM} and Yoshida's result in \cite{Yo} for an arbitrary Noetherian ring.

\begin{thm}
\label{homext} Let $R$ be a Noetherian ring and $I$ an ideal of $R$. Let $M$ be
a finitely generated $R$-module and $t\geq 1$ be a positive integer such that the
$R$-modules $H^i_I(M)$ are ${\rm FD_{\leq1}}$ ~ $R$-modules for all $i<t$.
Then, the following conditions hold:

{\rm(i)} The $R$-modules $H^i_I(M)$ are $I$-cofinite for all
$i<t$.

 {\rm(ii)} For all ${\rm FD_{\leq 0}}$ {\rm{(}}or minimax{\rm{)}} submodule $N$ of
$H^{t}_I(M)$, the $R$-modules $${\rm Hom}_R(R/I,H^{t}_I(M)/N)\,\,\, {\rm
and}\,\,\,{\rm Ext}^1_R(R/I,H^{t}_I(M)/N)$$ are finitely generated. In particular the
set $\Ass_R(H^{t}_I(M)/N)$ is a finite set.
\end{thm}
\proof (i) We proceed by induction on $ t$. By Lemma \ref{2.2} the case $ t = 1$
is obvious since $H^{0}_I(M)$ is finitely generated. So, let $ t >1 $ and the result
has been proved for smaller values of $t$. By the inductive assumption,
$H^{i}_I(M)$ is I-cofinite for $ i = 0, 1, ..., t-2$. Hence by Lemma \ref{2.2}  and
assumption, ${\rm Hom}_R(R/I,H^{t-1}_I(M))\,\,\, {\rm and}\,\,\,{\rm
Ext}^1_R(R/I,H^{t-1}_I(M))$ are finitely generated. Therefore by Corollary \ref{FD1} ,
$H^{i}_I(M)$ is I-cofinite for all $ i < t$. This completes the inductive step.\\
%Use induction on $i$ and corollary 2.3
%and Lemma 2.4.
(ii) In view of   (i) and lemma \ref{2.2},  ${\rm Hom}_R(R/I,H^{t}_I(M))\,\,\, {\rm
and}\,\,\,{\rm Ext}^1_R(R/I,H^{t}_I(M))$ are finitely generated. On the other hand,
according to Lemma \ref{FD0} or  Melkersson's result  \cite[Proposition
4.3]{Mel}, N is $I$-cofinite. Now, the exact sequence

$$0\longrightarrow N\longrightarrow  H^{t}_I(M) \longrightarrow H^{t}_I(M)/N
\longrightarrow 0$$

 induces the following exact sequence,
$${\rm
Hom}_R(R/I,H^{t}_I(M))\longrightarrow {\rm
Hom}_R(R/I,H^{t}_I(M)/N)\longrightarrow {\rm Ext}^1_R(R/I,N)\longrightarrow$$ $$
{\rm Ext}^1_R(R/I,H^{t}_I(M))\longrightarrow {\rm
Ext}^1_R(R/I,H^{t}_I(M)/N)\longrightarrow {\rm Ext}^2_R(R/I,N).$$

Consequently  $${\rm Hom}_R(R/I,H^{t}_I(M)/N)\,\,\, {\rm and}\,\,\,{\rm
Ext}^1_R(R/I,H^{t}_I(M)/N)$$ are finitely generated, as required.\qed\\

\begin{cor}
\label{cof} Let $R$ be a Noetherian ring and $I$ an ideal of $R$. Let $M$ be a
finitely generated $R$-module  such that the $R$-modules $H^i_I(M)$ are ${\rm
FD_{\leq1}}$ {\rm{(}}or weakly Laskerian{\rm{)}}~ $R$-modules for all $i$. Then,
the

{\rm(i)} The $R$-modules $H^i_I(M)$ are $I$-cofinite for all $i$.

 {\rm(ii)} For any $i\geq 0$ and for any ${\rm FD_{\leq 0}}$ {\rm{(}}or minimax{\rm{)}} submodule $N$ of
$H^{i}_I(M)$, the $R$-module $H^{i}_I(M)/N$ is  $I$-cofinite.
\end{cor}
\proof (i)  Clear.\\ (ii)  In view of  (i) the R-module $H^i_I(M)$ is $I$-cofinite for all
$i$. Hence the R-module ${\rm Hom}_R(R/I,N)$ is finitely generated, and so it
follows from Lemma \ref{FD0} or \cite[Proposition 4.3]{Mel} that N is $I$-cofinite .
Now, the exact sequence

$$0\longrightarrow N\longrightarrow  H^{t}_I(M) \longrightarrow H^{t}_I(M)/N
\longrightarrow 0$$,

implies that the R-module $H^{i}_I(M)/N$ is  $I$-cofinite.\qed\\

\begin{cor}
\label{wl}
Let $R$ be a Noetherian ring and $I$ an ideal of $R$. Let $M$ be a
finitely generated $R$-module and $t\geq 1$ be a positive integer such that the
$R$-modules $H^i_I(M)$ are weakly Laskerian for all $i<t$. Then, the following
conditions hold:

{\rm(i)} The $R$-modules $H^i_I(M)$ are $I$-cofinite for all $i<t$.

 {\rm(ii)} For all ${\rm FD_{\leq 0}}$ {\rm{(}}or minimax{\rm{)}} submodule $N$ of
$H^{t}_I(M)$, the $R$-modules $${\rm Hom}_R(R/I,H^{t}_I(M)/N)\,\,\, {\rm
and}\,\,\,{\rm Ext}^1_R(R/I,H^{t}_I(M)/N)$$ are finitely generated.  In particular the
set $\Ass_R(H^{t}_I(M)/N)$ is a finite set.
\end{cor}
\proof Use Theorem \ref{fsf} and note that the category of weakly Laskerian
modules is contained in the category of ${\rm FD_{\leq1}}$ modules.\qed\\

One of the main result of this section is to prove that for an arbitrary ideal $I$ of a
Noetherian ring $R$, the Category
 of $I$-cofinite ${\rm FD_{\leq1}}$ modules compose an Abelian category.

\begin{thm}
\label{abel}
 Let $I$ be an ideal of a Noetherian ring $R$. Let $\mathscr {FD}^1(R,I)_{cof}$ denote the
category of $I$-cofinite ${\rm FD_{\leq1}}$ ~ $R$-modules. Then  $\mathscr
{FD}^1(R,I)_{cof}$  is an Abelian category.
\end{thm}
\proof Let $M,N\in \mathscr {FD}^1(R,I)_{cof}$ and let $f:M\longrightarrow N$ be
an $R$-homomorphism. It is enough that to show that the $R$-modules $\ker f$
and ${\rm coker} f$ are $I$-cofinite.

To this end, the exact sequence
$$0\longrightarrow \ker f \longrightarrow M \longrightarrow
{\rm im} f \longrightarrow 0,$$ induces an  exact sequence
$$0\longrightarrow {\rm Hom}_R(R/I,\ker f) \longrightarrow {\rm
Hom}_R(R/I,M) \longrightarrow {\rm Hom}_R(R/I,{\rm im} f)$$$$\longrightarrow
{\rm Ext}^{1}_R(R/I,\ker f) \longrightarrow \Ext^{1}_R(R/I,M),$$ that implies the
$R$-modules $\Hom_R(R/I,\ker f)$ and $\Ext^1_R(R/I,\ker f)$ are finitely
generated. Therefore it follows from Theorem \ref{FD1} that $\ker f$ is $I$-cofinite.
Now, the assertion follows from the following exact sequences
$$0\longrightarrow \ker f \longrightarrow M \longrightarrow
{\rm im}f \longrightarrow 0,$$and
$$0\longrightarrow {\rm im}f \longrightarrow N \longrightarrow
{\rm coker}f \longrightarrow 0.$$\qed\\

The following corrollary is a generalization of \cite[Theorem 2.7]{AB}.

\begin{cor}
\label{torext}
 Let $I$ be an ideal of a Noetherian ring $R$. Let  $M$ be an ${\rm FD_{\leq1}}$ ~ $I$-cofinite  $R$-module.
  Then, the R-modules ${\rm Ext}^{i}_R(N,M)$ and ${\rm Tor}^R_{i}(N,M)$
are $I$-cofinite and ${\rm FD_{\leq1}}$ modules, for all finitely generated
R-modules $N$ and all integers $ i \geq 0$.
\end{cor}
\proof Since $N$ is finitely generated it follows that $N$ has a free resolution of
finitely generated free modules. Now the assertion follows using Theorem
\ref{abel} and computing the modules ${\rm Tor}_i^R(N,M)$ and ${\rm
Ext}^i_R(N,M)$, by this
free resolution. \qed\\

\begin{cor}
\label{torextlc} Let $I$ be an ideal of  a Noetherian ring R, M a non-zero finite
$R$-module such that $\dim M/IM\leq 1$\rm (e.g., $ \dim R/I\leq 1)$. Then for
each finite $R$-module $N$, the $R$-modules ${\rm Ext}^{j}_R(N,H^{i}_I(M))$ and
${\rm Tor}^R_{j}(N,H^{i}_I(M))$ are $I$-cofinite for all $ i \geq 0$ and $ j \geq 0$.
\end{cor}
\proof Note that $\dim \Supp H^{i}_I(M) \leq \dim M/IM \leq 1$ thus it is an ${\rm
FD_{\leq1}}$ module for all  $ i \geq 0$, now use  Corollary \ref{torext}.
\qed\\

\begin{lem}
Let $R$ be a Noetherian ring, $I$ a proper ideal of $R$ and $M$ be a non-zero
 ${\rm D_{\leq1}}$ and $I$-cofinite  $R$-module. Then for each non-zero finitely
generated $R$-module $N$ with support in $V(I)$, the $R$-modules
$\Ext^i_R(M,N)$ are finitely generated, for all integers $i\geq 0$.
\end{lem}
\proof See \cite[Theorem 2.8]{IGB}.\qed\\

\begin{cor}
 Let $R$ be a Noetherian ring and $I$ be an ideal of  $R$.  Let  $M$ be an ${\rm FD_{\leq1}}$ and  $I$-cofinite  $R$-module.
  Then, the R-modules ${\rm Ext}^{i}_R(M,N)$ and ${\rm Tor}^R_{i}(M,N)$
are finitely generated, for all finitely generated R-modules $N$ with ${\rm
Supp}(N)\subseteq V(I)$ and all integers $ i \geq 0$.
\end{cor}
\proof The assertion follows from the definition using Lemma 3.9 and \cite[Theorem 2.1]{Mel}.\qed\\

\end{document}